\documentclass{article}    
\usepackage{amstext}
\def\qed{$\rlap{$\sqcap$}\sqcup$}
\usepackage{latexsym}

\begin{document}           

{\ }\\
\begin{center}
{\huge {\bf Stanley's theorem on codimension 3\\Gorenstein $h$-vectors}} \\ [.250in]
{\large FABRIZIO ZANELLO\\
Dipartimento di Matematica, Universit\`a di Genova, Genova, Italy.\\E-mail: zanello@dima.unige.it}
\end{center}

{\ }\\
\\
ABSTRACT. In this note we supply an elementary proof of the following well-known theorem of R. Stanley: the $h$-vectors of Gorenstein algebras of codimension 3 are SI-sequences, i.e. are symmetric and the first difference of their first half is an $O$-sequence.

{\large

{\ }\\
\\\indent
We consider standard graded artinian algebras $A=R/I$, where $R=k[x_1,...,x_r]$, $k$ is any field, the $x_i$'s have degree 1 and $I$ is a homogeneous ideal of $R$. Recall that the {\it $h$-vector} of $A$ is $h(A)=h=(h_0,h_1,...,h_e)$, where $h_i=\dim_k A_i$ and $e$ is the last index such that $\dim_k A_e>0$. Since we may suppose that $I$ does not contain non-zero forms of degree 1, $r=h_1$ is defined to be the {\it codimension} of $A$.\\\indent 
The {\it socle} of $A$ is the annihilator of the maximal homogeneous ideal $\overline{m}=(\overline{x_1},...,\overline{x_r})\subset A$, namely $soc(A)=\lbrace a\in A {\ } \mid {\ } a\overline{m}=0\rbrace $. Since $soc(A)$ is a homogeneous ideal, we define the {\it socle-vector} of $A$ as $s(A)=s=(s_0,s_1,...,s_e)$, where $s_i=\dim_k soc(A)_i$. Note that $s_e=h_e>0$. The integer $e$ is called the {\it socle degree} of $A$ (or of $h(A)$). If $s=(0,0,...,0,s_e=1)$, we say that the algebra $A$ is {\it Gorenstein}.\\\indent
The next theorem is a well-known result of Macaulay.\\
\\\indent 
{\bf Definition-Remark 1.} Let $n$ and $i$ be positive integers. The {\it i-binomial expansion of n} is $$n_{(i)}={n_i\choose i}+{n_{i-1}\choose i-1}+...+{n_j\choose j},$$ where $n_i>n_{i-1}>...>n_j\geq j\geq 1$.\\\indent 
Under these hypotheses, the $i$-binomial expansion of $n$ is unique (e.g., see $[BH]$, Lemma 4.2.6).\\\indent 
Furthermore, define $$n^{<i>}={n_i+1\choose i+1}+{n_{i-1}+1\choose i-1+1}+...+{n_j+1\choose j+1}.$$
\\\indent 
{\bf Theorem 2} (Macaulay). {\it Let $h=(h_i)_{i\geq 0}$ be a sequence of non-negative integers, such that $h_0=1$, $h_1=r$ and $h_i=0$ for $i>e$. Then $h$ is the $h$-vector of some standard graded artinian algebra if and only if, for every $d$, $1\leq d\leq e-1$, $$h_{d+1}\leq h_d^{<d>}.$$}
\\\indent 
{\bf Proof.} See $[BH]$, Theorem 4.2.10. (This theorem holds, with appropriate modifications, for any standard graded algebra, not necessarily artinian.){\ }{\ }\qed \\
\\\indent
Let us now recall a few definitions. We denote by $\lfloor x\rfloor $, as usual, the largest integer less than or equal to $x$.\\
\\\indent
{\bf Definition 3.} i). A sequence of non-negative integers which satisfies the growth condition of Macaulay's theorem is called an {\it $O$-sequence}.\\\indent
ii). A vector of non-negative integers $v=(v_0,v_1,...,v_d)$ is {\it differentiable} if its first difference, $$\Delta v=((\Delta v)_0=1,(\Delta v)_1=v_1-v_0,...,(\Delta v)_d=v_d-v_{d-1}),$$ is an $O$-sequence. (It is easy to see that if $v$ is differentiable, then $v$ is itself an $O$-sequence.)\\\indent
iii). An $h$-vector $h=(1,h_1,...,h_e)$ is an {\it SI-sequence} (named after Stanley and Iarrobino), if it is symmetric with respect to ${e\over 2}$ and if its first half, $(1,h_1,...,h_{\left\lfloor {e\over 2}\right\rfloor })$, is differentiable.\\
\\\indent
The study of the possible Gorenstein $h$-vectors is a central problem in Commutative Algebra. Initially, Stanley and Iarrobino (independently) conjectured that all Gorenstein $h$-vectors (of any codimension $r$) are SI-sequences. (See also Harima's paper $[Ha]$ on the $h$-vectors of Gorenstein algebras having the weak Lefschetz property.)\\\indent
The fact that Gorenstein $h$-vectors are symmetric is well-known, and a theorem of Migliore-Nagel ($[MN]$) and Cho-Iarrobino ($[CI]$) shows that, in any codimension, an SI-sequence is always a Gorenstein $h$-vector.\\\indent
The converse (that, as we just said, was conjectured to be true), instead has been proven false for $r\geq 5$. In particular, not even all Gorenstein $h$-vectors are {\it unimodal} (i.e. they do not increase after they start decreasing). The first example of a Gorenstein algebra with a non-unimodal $h$-vector was given by Stanley (see $[St1]$, Example 4.3), in codimension 13. Later, Bernstein-Iarrobino ($[BI]$), Boij-Laksov ($[BL]$) and Boij ($[Bo]$) exhibited many other non-unimodal Gorenstein $h$-vectors of codimension 5 or greater.\\\indent
In codimension 4, we do not know whether or not all Gorenstein $h$-vectors $h$ are SI-sequences, nor even whether they must be all unimodal. There is, however, a remarkable result of Iarrobino and Srinivasan ($[IS]$) which shows that, if the entry of degree 2 of $h$ is less than or equal to 7, then $h$ must be an SI-sequence.\\\indent
Instead, in codimension 2, the conjecture that all Gorenstein $h$-vectors are SI-sequences is correct, as first observed by Macaulay ($[Ma]$), and is indeed an easy exercise assuming Macaulay's Theorem 2 and the symmetry of Gorenstein $h$-vectors.\\\indent
In codimension 3, the conjecture still holds true, as shown by Stanley (see $[St1]$, Theorem 4.2). His proof is based on a deep structure theorem due to Buchsbaum and Eisenbud ($[BE]$, Proposition 3.3). The purpose of the present note is to supply an elementary proof of this important result of Stanley.\\
\\\indent
{\bf Theorem 4.} {\it Let $h$ be a Gorenstein $h$-vector of codimension 3. Then $h$ is an SI-sequence.}\\
\\\indent
Before going into the proof we need to recall the following observation of Stanley:\\
\\\indent
{\bf Remark 5} (Stanley). {\it Let $h=(1,h_1,...,h_e)$ be a Gorenstein $h$-vector. Then, for any index $j\geq 1$, there exists a Gorenstein $h$-vector $(a_j=1,a_{j+1},...,a_e)$ such that the vector $(1,h_1,...,h_{j-1},h_j-a_j,...,h_e-a_e)$ is an $O$-sequence.}\\
\\\indent
{\bf Proof.} See $[St2]$, bottom of p. 67.{\ }{\ }\qed \\
\\\indent
{\bf Proof of Theorem 4.} Let $h=(1,h_1=3,h_2,...,h_e)$ be a Gorenstein $h$-vector of codimension 3. We want to show that its first half is differentiable. By Remark 5 (with $j=1$), there exists a Gorenstein $h$-vector $a=(a_1=1,a_2,...,a_{e-1},a_e)$ such that $$(1,\Delta_1=2,\Delta_2=h_2-a_2,...,\Delta_e=h_e-a_e)$$ is an $O$-sequence. Notice that, by this choice of the indices, $a$ is symmetric with respect to ${e+1\over 2}$; in particular, $a_2=a_{e-1}\leq h_{e-1}=h_1=3$.\\\indent
First we show that $h$ is unimodal. Suppose it is not. We may assume, by induction, that $e$ is the least socle degree for which there exists a non-unimodal Gorenstein $h$-vector of codimension 3. Hence we have $h_i<h_{i-1}$ for some $i\leq \left\lfloor {e\over 2}\right\rfloor $. Since $a$ is unimodal (by the induction hypothesis, if it has codimension 3), we have $\Delta_i<i+1$ (i.e. $\Delta_i$ is not {\it generic}). But $$\Delta_{e-(i-1)}=h_{e-(i-1)}-a_{e-(i-1)}=h_{i-1}-a_i>h_i-a_i=\Delta_i,$$ a contradiction, since, by Macaulay's Theorem 2, an $O$-sequence starting with $(1,2,...)$ cannot increase after it is no longer generic. This proves that $h$ is unimodal.\\\indent
Now we want to show that the first half of $h$ is differentiable. We may suppose, by induction, that all Gorenstein $h$-vectors of codimension 3 and socle degree lower than $e$ are SI-sequences. The differentiability of $h$ is obvious as long as $h$ is generic (i.e. $h_i={i+2\choose 2}$). Hence suppose that $h_i$ is not generic, for some $i\leq \left\lfloor {e\over 2}\right\rfloor $. By Macaulay's theorem, we need to show that \begin{equation}\label{1}h_i-h_{i-1}\leq h_{i-1}-h_{i-2}.\end{equation}\indent
Let us first consider the case $a_i={i+1\choose 2}$. We have $${i+1\choose 2}=a_i=a_{e-(i-1)}\leq h_{e-(i-1)}=h_{i-1}\leq {i+1\choose 2},$$ and therefore $h_{i-1}={i+1\choose 2}$, i.e. $h_{i-1}$ is generic. Thus, (\ref{1}) becomes $h_i-{i+1\choose 2}\leq i$, and this is true since $h_i<{i+2\choose 2}$. This proves the theorem for $a_i={i+1\choose 2}$.\\\indent
Hence, let us assume from now on that $a_i<{i+1\choose 2}$. Suppose now that $\Delta_{i-1}$ is generic (i.e. $\Delta_{i-1}=i$). Therefore (\ref{1}) becomes $$a_i+\Delta_i-a_{i-1}-i\leq a_{i-1}+i-a_{i-2}-(i-1),$$ which is true since $\Delta_i-i\leq 1$ and $a_i-a_{i-1}\leq a_{i-1}-a_{i-2}$, because $a_i<{i+1\choose 2}$ and $a$ is an SI-sequence (by induction, if it has codimension 3). This completes the proof for $\Delta_{i-1}=i$.\\\indent
Hence let us suppose that $\Delta_{i-1}\leq i-1$. Therefore $\Delta_{i-1}\geq \Delta_{e-(i-1)}\geq \Delta_{e-(i-2)}$, whence $$a_{e-(i-1)}-a_{e-(i-2)}\leq h_{e-(i-1)}-h_{e-(i-2)},$$ i.e. \begin{equation}\label{2}a_i-a_{i-1}\leq h_{i-1}-h_{i-2}.\end{equation}\indent
Similarly, $\Delta_{i-1}\geq \Delta_i$, i.e. \begin{equation}\label{3}h_i-h_{i-1}\leq a_i-a_{i-1}.\end{equation}\indent Thus, (\ref{1}) follows from (\ref{3}) and (\ref{2}). This proves the theorem.{\ }{\ }\qed \\

\clearpage

{\bf \huge References}\\
\\
$[BI]$ {\ } D. Bernstein and A. Iarrobino: {\it A nonunimodal graded Gorenstein Artin algebra in codimension five}, Comm. in Algebra 20 (1992), No. 8, 2323-2336.\\
$[Bo]$ {\ } M. Boij: {\it Graded Gorenstein Artin algebras whose Hilbert functions have a large number of valleys}, Comm. in Algebra 23 (1995), No. 1, 97-103.\\
$[BL]$ {\ } M. Boij and D. Laksov: {\it Nonunimodality of graded Gorenstein Artin algebras}, Proc. Amer. Math. Soc. 120 (1994), 1083-1092.\\
$[BH]$ {\ } W. Bruns and J. Herzog: {\it Cohen-Macaulay rings}, Cambridge studies in advanced mathematics, No. 39, Revised edition (1998), Cambridge, U.K..\\
$[BE]$ {\ } D.A. Buchsbaum and D. Eisenbud: {\it Algebra structures for finite free resolutions, and some structure theorems for ideals of codimension 3}, Amer. J. Math. 99 (1977), 447-485.\\
$[CI]$ {\ } Y.H. Cho and A. Iarrobino: {\it Inverse Systems of Zero-dimensional Schemes in ${\bf P}^n$}, J. of Algebra, to appear.\\
$[Ha]$ {\ } T. Harima: {\it Characterization of Hilbert functions of Gorenstein Artin algebras with the weak Stanley property}, Proc. Amer. Math. Soc. 123 (1995), No. 12, 3631-3638.\\
$[IS]$ {\ } A. Iarrobino and H. Srinivasan: {\it Some Gorenstein Artin algebras of embedding dimension four, I: components of $PGOR(H)$ for $H=(1,4,7,...,1)$}, J. of Pure and Applied Algebra, to appear.\\
$[Ma]$ {\ } F.H.S. Macaulay: {\it The Algebraic Theory of Modular Systems}, Cambridge Univ. Press, Cambridge, U.K. (1916).\\
$[MN]$ {\ } J. Migliore and U. Nagel: {\it Reduced arithmetically Gorenstein schemes and simplicial polytopes with maximal Betti numbers}, Adv. in Math. 180 (2003), 1-63.\\
$[St1]$ {\ } R. Stanley: {\it Hilbert functions of graded algebras}, Adv. Math. 28 (1978), 57-83.\\
$[St2]$ {\ } R. Stanley: {\it Combinatorics and Commutative Algebra}, Second Ed., Progress in Mathematics 41 (1996), Birkh\"auser, Boston, MA.

}

\end{document}